# Gabriel Quotient Rings

**By C.L.Wangneo**

**Jammu ( J-K-UT)**

**India -180002 .**

**Abstract :--**In this paper we construct and initiate the study of the Gabriel Quotient ring of a prime noetherian ring R at a m-gabriel filter g . The statement of the theorem that we prove is given below ;

**Theorem ( construction and definition ) :-** Let R be a prime noetherian ring with k.dim. (R) = n , n a positive integer. Let Q denote the Goldie quotient ring of R . For a fixed positive integer m, 0 <m < n , call xm = {all pεspec(R) / k-dim. (R/p))=m} a full set of prime ideals of R . Let cm = { all cε R / k-dim. (R/cR )< m } , and let g = {all right ideals J ≤ R / k-dim. (R/ J )< m } be the m-Gabriel filter in R . We then construct an extension ring R(m) of R at the m-Gabriel filter g with the following properties (i) R(m) is a right noetherian subring of Q with identity element 1 of R . (ii) If u(R(m)) is the set of units of R(m) then u(R(m) п R = cm precisely.

(iii) For a full set of m-prime ideals the set cm is a right ore set .

**Introduction :-** This paper is divided into three sections . In section (1) we introduce glossary of terms of use to us throughout this paper . This glossary of terms appears already as preliminaries in [3] but here we present it briefly for the sake of completeness.Then in section (2) we mention further glossary of terms and results that lead to the definition and construction of the Gabriel extension ring of a prime noetherian Ring R . Then in section (3) we prove that the Gabriel extension ring of section (2) is actually a quotient ring of the ring R . This paper should be considered as a brief replacement of [4] .

**<u>Notation and Terminology</u> :**

**(a) Throughout in this paper a ring is meant to be an associative ring with identity which is not necessarily a commutative ring.**

**(b)We will throughout adhere to the same notation and terminology for the definitions occurring in this paper as in [3].**

**By a module M over a ring R we mean that M is a right R-module unless stated otherwise. For the basic definitions regarding right noetherian modules over right noetherian rings and the definitions of uniform modules , right Krull dimension , critical modules , krull homogenous modules we refer the reader to both [1] and [2 ]. Also one may consult both [1] and [2[ for the definition and properties of Ass.(M) and its associated concepts such as faithful and fully faithful module over a right noetherian ring Moreover we will use the following terminology throughout .**

**(c) If R is a ring then we denote by Spec.R , the set of prime ideals of R . Also right Krull dimension of the right R-moodule M if it exists will throughout be denoted by |M| . For two subsets A and B of a given set , A ≤ B means B contains A and A < B denotes A ≤ B but A≠B. Also for two sets A and B , A⊄B denotes the set B that does not contain the subset A. For an ideal A of R , c(A) denotes the set of elements of R that are regular modulo A . Recall if S ≤ M is a non-empty subset , then we denote the right annihilator of S in R by r (S) . Also recall that if R is a ring and M is a right R module and if T is a multiplicatively closed subset of regular elements of R then a submodule N of M is said to be T tosion if for any element x of M there exists an element t in T such that x t = 0 . M is said to be T torsion free if for any nonzero element x of M and for any t in T , x t ≠ 0 . If M is not T- torsion , then there exists at least one non-zero element x of M such that x is a T Torsion free element of M . This means x t ≠ 0 for some element t ≠ 0 of T . Similarly if M is not T-**

**torsionfree , then there exists at least one non-zero element y of M such that y is a T- Torsion element of M . This means yt = 0 , for at least one element t≠ 0 in T . In this case then we say that M is a T-non-torsion right R module .**

**We will mention a few words below about the standard definition of the extension ring of a ring R to another ring S containing it and for the definition of the extension of a right ideal A of R to S and the contraction of a right ideal B of S to R .**

**Let R be a prime noetherian ring and let Q denote its goldie quotient ring . Then we will call S an extension ring of R in Q if (i) S is a subring of Q such that R≤ S ≤Q (ii) $1_R = 1_s$ is the identity element of S . The standard notation for an extension ring S of R is S= $R^e$ Moreover then for any right ideal A of R the extension of the right ideal A to S is denoted by $A^e$ which is defined as the right ideal of S such that $A^e$ = AS , that is $A^e$ is the right ideal generated by A in S . Similarly for any right ideal B of S we write $B^C = B \cap R$ . Clearly $B^C$ is a right ideal of R . We call $B^c$ the contraction of the right ideal B of S=$R^e$ to R**

**We will throughout write Spec.(R) = family of prime ideals of R , max. specR = family of maximal ideals of R and prim spec. R =family of right primitive ideals of R .**

**A word about the numbering of a theorem , propostion ,lemma etc. All such numbering occurs on the left of each proposition , lemma etc. but while referring to these theorems lemmas etc. we refer them in the usual way in the paper For example (1.2) Lemma would be referred as Lemma (1.2) in the paper while quoting it .in the paper .**

**As mentioned in the introduction this paper has three sections namely , section (1) , section (2) and section (3) . In section (1) we mention the glossary of terms of use to us always as they appear in section (1) of [3] . In section (2) we define and construct an**

extension ring of a prime noetherian ring R at a m-gabriel filter g . Then in section (3) we complete our main theorem of this paper.

**(1) section ( Glossay of terms ) :-** In this section we first revise briefly the definitions , results and symbols of use to us in section (3) in which we prove our main theorem that introduces and constructs the Gabriel quotient ring of a prime noetherian ring R at a m-gabriel filter . From section (1) of [3] especially from (1.1) to (1.3) , we have the following ;

**(1.1) If** R is a prime , noetherian ring with ]R] = n , $n > 0$ then for a fixed positive integer m , $0 < m < n$ , we set xm = { All p ε spec. R / |R/ P| = m } . We call xm as the full set of m-prime ideals of R. Then for each m and its corresponding Xm we let vm=∩ c(p) ; for each p ε xm. and let cm = { all c ε R / IR/cR I< m } Then cm is a multiplicatively closed subset of regular elements of R and cm ≤ vm and vm is also a multiplicatively closed subset of regular elements of R

**(1.2)** For the prime noetherian ring R with IR I= n , $n > 0$ and for a positive integer m , $m<n$ , let xm , vm and cm be defined as in definition (a) above . We then have the following families of right ideals of R , namely ;

c = { Right ideals I of R / I∩Cm ≠ϕ} .

v = { Right ideals I of R / I∩ Vm ≠ϕ} .

g = { Right ideals I of R / IR/I I< m }.

h = { Right ideals I of R / I∩ c(p) ≠ϕ for each pεXm}

**(1.3)** For the prime noetherian ring R with IR I=n , $n > 0$ if xm , vm,cm and v , c , g and h are as in definitions (1.1 and (1.2)above then we say xm has the right intersection condition if for any right ideal I of R , $I \cap c(p) \neq \phi$ for each $p \varepsilon X_m$, implies that $I \cap V_m \neq \phi$. Note this is equivalent to saying that V=h, because $v \leq h$ always

(1.4) From lemma (1.2) of [3] we have that if ]R] =n, n, a finite integer , n > 0 then for a fixed positive integer . m . 0 < m < n let $x_m$, $v_m$ , $c_m$ be as in the above definition (1.1) . Now consider the families of right ideals c, v, g and h as in definition (1.2) above . Then the following are true ;

(i) $c_m \leq v_m$ and (ii) $C \leq g \leq h$ . Also $c \leq v \leq h$.

(iii) If $X_m$ has the right intersection condition(see definition (1.3) above ) then we have

(x) v= h and

(y)$g \leq v$ .

(1.5) (Remark) For notational convenience we will denote the Gabriel extension ring of a noetherian prime ring R at a m-gabriel filter g by R(m) . All the above definitions and notation including definitions (1.1) and (1.2) and the notation therein as well as the intersection condition remain valid if we replace the full set of prime ideals $X_m$ by any subset $Y_m$ of $X_m$ .In connection with this we would like to know if a full set of m-prime ideals $x_m$ of R has the desired property that the set $c_m$ associated to $X_m$ is a right ore set .

(1.6) Gabriel filters :- We now define m-Gabriel filters and some allied concepts ; . -

(a) Definition :- Let R be a ring . Let $g = \{ I < R / I$ is a right ideal $\}$ . Let g be a non-empty family . Then g is called a m-Gabriel filter if g has the following properties ;

(i) For any $I \varepsilon g$ if J is a right ideal such that $I \leq J$ , then $J \,\mathcal{E}\, g$ .

(ii) For I, J in g , I ∩ J Ɛ g .

(iii) For any a in R , and I ε g , the right ideal a-1 I = { x in R / ax in I } , is also in g .

Note if J = a-1 I , then J is a maximal right ideal of R such that a J < I .

(b) Definition In a ring R a multiplicative Gabriel filter is a nonempty family g of right ideals of R such that

(i) g is a Gabriel filter

(ii) For I, J in g the product of the right ideals IJ Ɛ g .

(1.7) We stick to the definitions and allied concepts of m-Gabriel filter from definition (1.4) to (1.6) of [3] . Thus if g is an m- Gabriel filter and If cm = {c in R / IR./ cRI < m } and c= { all right ideals k in R / K п cm ≠ ф } then given cm , c, and g we see that g=c implies cm is a right ore-set always .

(1.8) Similarly if xm is a full set of m-prime ideals of R and if vm = п c(p) for all p ε xm we denote by v= { all right ideals k in R / K п vm ≠ ф } then given vm , v, and g we see that g =v implies that vm is a right ore set and vm=cm and hence g=c .

(1.9) ( Some facts about the ring Rc ) :- We now .state the following facts about the ring Rc

(a) Let g=c . Clearly then cm is aright ore set . If Rc denotes the quotient ring of R at the right ore set cm then we have the following as true ;

(i) Given Q, R ,c and g=c as defined above then the extension ring Rc is a subset of Q such that Rc = { q ε Q / qc ε R / for some c ε cm).

Hence Rc is an extension ring of R in Q . Moreover then Jc =Rc iff J ε g where Jc is the usual extension of the right ideal J of R in Rc .

**(ii) If J is a right ideal of R , then Jc is a proper right ideal of Rc iff J п c=ф .**

**(iii) Hence Jc is a proper right ideal of Rc iff IR/JI ≥ m .**

**(iv) From (iii) above we get that if J is a proper right ideal of R such that R/J is a critical right R module then IR/JI = m and hence R/ J is a m-critical right R –module.**

**(v) From (iv) above we have that that Jc is a maximal right ideal of Rc ( or equivalently Rc/Jc is a simple right Rc module ) iff R/J is a m-critical right R-module . In particular we have that r-ann. (Rc/Jc)=qc iff r-ann. R/J =q . Hence qc and q are prime ideals of Rc and R respectively .**

**(vi) From (v) we infer that Rc/ Jc is either a trorsion free right Rc/ qc – module in which case then Rc /qc is a simple artinian ring or Rc/ Jc is a torsion right Rc/ qc – module in which case then Rc /qc can never be an artinian ring .**

**(vii) Moreover in (vi) above we have that if R/J is a torsion free R/q module then IR /.J I= m =IR/qI or else if R/J is a torsion R/q module then we must have that IR /.J I= m < IR/qI.**

**(1.10)(a ) Given Q, R ,c and g=c as defined in (1.9) above then recall the extension ring Rc is a subset of Q such that Rc = { q ε Q / qc ε R / for some c ε cm). If we call an element q'=q+R of the right R module Q /R as a cm or tc torsion element of Q/R provided qc ε R , for some c ε cm . then we see that Jc is a proper right ideal of Rc iff J п c=ф**

**(b) ( Theorem ) :- With the notaion in (a) above given Q, R,c and g=c as above then we have the following two descriptions of Rc given as below ;**

**(i) Rc is a R- submodule of Q such that R ≤ Rc and and Rc = { q ε Rc / qc ε R for some c ε cm }**

(ii) if q' = q+R denotes the image of q in Rc/R then the submodule q'R of Rc/R is such that Iq'RI < m . Hence Rc/R ={ all q ε R / if q' = q+R then Iq'RI < m } .Thus in this case each q ε Rc can be recovered from the preimage q of q' in Rc/R } .

**Hence we have the following definition ;**

(c) (Definition ) the above theorem yields the following definition of two equivalent descriptions of the right R-module Rc given below

(i) Rc is a R- submodule of Q such that R ≤ Rc and and Rc = { q ε Rc / qc ε R for some c in ε cm }

(ii) Rc/R = tc(Q/R ) = {all elements q in Q / the image q'=q+R of q in Q/R has right krull dimension , Iq'RI < m }.

**(2) Section (Gabriel extension ring ):---** In this section we construct an extension ring R(m) of a prime noetherian ring R with some useful properties .

**(2.1) Definition(Gabriel extension ring ):---** We follow the description of Rc in (1.10) above to describe a general torsion theoretic extension ring $R^e$ of a noetherian prime ring at a m-gabriel filter g .We call this exrension ring of R as the gabriel extension ring . We construct the gabriel extension ring of R in three steps given below .

(i) (step (1) ) ;- Let R be a prime noetherian ring with I R I = n, n , a positive integer and let Q be its goldie quotient ring . Let c(0) be the set of regular elements of R . For a fixed positive integer m , m < n , let g be the m- Gabriel filter defined as g = { all right ideals I in R / | R/I | < m } .we call an extension ring $R^e$ of R in Q a Gabriel extension ring of R if the following hold true ;

(i) $R^e$ = { q ε Q / qJ ≤ R / for some J ε g ) . It is clear then that $R^e$ is an R- submodule of Q such that R ≤ $R^e$ and $R^e$ /R = tg(Q/R) = tg-torsion R submodule of Q/R . Note that for any element q'= q+R , q in Q is called a tg-torsion element if qJ ε R for some J in g . Moreover then Q / $R^e$ is a tg-torsion-free factor module of Q/R .

(ii) ( step (2) :- In this step we describe an alternative description of the set $R^e$ . For each q in Q , if q ' = q +R is image of the element q in Q/R and if q' R = right submodule of Q/ R generated by q' in Q /R then an alternative description of $R^e$ / R is that $R^e$ / R is a right R submodule of Q /R such that $R^e$ / R = { q in Q/|q'R |< where q'=q+R is the image of the element q in Q/R } . In this case then $R^e$ is the preimage of $R^e$ /R in Q . Moreover then for any nonzero element q in Q if q' = q+ $R^e$ is nonzero , then Iq'RI ≥ m .

(iii) (step(3) :- It is not difficult to see from steps (1 and (2) above that $R^e$ = { q in Q / qJ < R , for some J in g } is a subring of Q with addition and multiplication inherited from Q naturally .We prove this as a theorem below .

(2.2) Theorem (construction and definition ):- (a) The Gabriel extension module $R^e$ of a prime noetherian ring R at the m-Gabriel filter g is defined and denoted as the right R- sub-module of Q , such that . $R^e$ = { q in Q / qJ < R , for some J in g } which can also be recovered from the description of $R^e$ / R = { q in Q/|q'R |< m, where q'=q+R in Q/R } by considering each element of $R^e$ as the preimage of elements in Q from $R^e$ / R . We now show that $R^e$ is a subring of Q with addition and multiplication inherited from Q naturally .

Proof :- Note from steps (1) and (2) above we have that for each q in $R^e$ we have that q J ≤ R . From step (2) it follows that if q ' = q +R is image of the element q in Q/R and if q' R = right submodule of Q/ R generated by q' in Q /R then from step (1) an alternative description of the extension ring $R^e$ is that $R^e$ = { q in Q/|q'R |< m, where q' = q+R} . We now show if $R^e$ = { q in Q / qJ < R , for some J in g } then $R^e$ is a subring of Q . Denote for any element q

in Q its image in the right R module Q/R by q’ . Now we show that for any two elements q and t in $R^e$ if their images in Q /R are q’,t’ ε $R^e$ . Let J , H in g be such that q J ≤ R and t H ≤ R . So by the description of $R^e$ above we have that Iq’RI <m and It’RI <m . Now qtR ≤ qR as right R modules . Hence Iq’t’RI ≤ Iq’RI< m , where q’t’R and q’R are considered as right R sub-modules of Q/R . Hence using this second description of elements of $R^e$ above we have since q’t’R in Q / R is such that Iq’t’RI < m hence the preimage of q’t’ in Q/R is qt Thus qt must belong to $R^e$ . Hence $R^e$ is a ring .

**(2.3) Thus we define and call** $R^e$ the Gabriel extension ring of R at the m-gabriel filter g if $R^e$ = { q ε Q / q J ≤ R , for some J ε g }. Moreover another description of $R^e$ is that for each q in Q , if q ‘ = q +R is image of the element q in Q/R and if q’ R = right submodule of Q/ R generated by q’=q+R in Q /R then $R^e$ = { q ε Q which are preimages of the elements q’ =q+R in Q/R with |q’R |< m} . .

**(2.4) Theorem ( Units of the Gabriel extension ring) :** -In the theorem below we describe the units of the Gabriel quotient ring at the filter g as described above .

Let R be a prime , noetherian ring with I R I = n, n , a finite integer , n > 0 and let Q be its Goldie extension ring . For a fixed integer , m, m <n , let g be the m- Gabriel filter defined as g = { all right ideals I in R / | R/I | < m } . Let $R^e$ be the Gabriel extension ring of R at the Gabriel filter g . Set now cm ={elements c in R / | R/cR | < m } . Let U($R^e$)=Set of all units of $R^e$ . Then we have that cm is precisely the set cm = U($R^e$) ∩ R .

**Proof :**- The proof is on the following lines . Let u ε $R^e$ π R where u is a unit of $R^e$ . Then $u^{-1}$ ε $R^e$ implies that there is a J ε g such that $u^{-1}$ J ≤ R .Since $R^e$ is a ring , so multiplying on the left by u we get that J ≤ u R . Since J ε g , hence u R ε g as well . Thus by the definition of cm we get that u ε cm . This proves that U($R^e$) ∩ R ≤ cm We now show that cm ≤ U($R^e$) ∩ R. But this is obvious because the goldie rank of $R^e$ is finite . Hence for a c ε c(0) implies since $c^{-1}$ ε Q and cR ≤ R thus $c^{-1}$ ε $R^e$ is the left inverse of c in cm . But then it is obvious that $c^{-1}$ is also the right inverse of c . Hence c is a unit of $R^e$ as well . Thus U($R^e$) ∩ R.≤ cm . Hence this proves that U($R^e$) ∩ R = cm .

**(3)section** :- In this section we show that the set cm is indeed a right ore set .

**(3.1) Remark** :- For notational convenience we will denote the Gabriel quotient ring of a noetherian prime ring R at a m-gabriel filter g by R(m) . We see that if U(R(m) is the set of units of R(m) then U(R(m) п R = cm precisely . Obviously then one can ask as to when is cm a right ore set . In connection with this we would like to consider the full set of m prime ideals namely Xm thus enlaging the set of prime ideals for which cm is regular module each p in Xm . We will see that it tuns out to be useful .

**(3.2) Definition :-** Let R be a prime , noetherian ring with I R I = n, n , a positive integer and let Q be its Goldie quotient ring . Let c(0) be the -set of regular elements of R . For a fixed positive integer , m, m < n , let g be the m- Gabriel filter defined as g = { all right ideals I in R / | R/I | < m} . Let cm = { all c ε R / | R/cRI | < m} and let c = { family of right ideals k of R / k п Cm ≠ ϕ } . Let xm = { all prime ideals p ε R / | R/p | = m } be the full set of m-prime ideals of R . Let for each p ε xm , c(p) have the usual meaning . We denote for each p ε xm the family gp = { family of right ideals J in g / J п c(p) ≠ ϕ } . Let $R^e$ = { q in Q / qJ < R , for some J in g } be the Gabriel extension ring of R at the Gabriel filter g. Then we define the following ;

We say the filter g is a localising filter for a prime ideal p of R if g ≤ g(p) .

**(3.3) Theorem :-** Let R be a prime , noetherian ring with I R I = n, n , a positive integer and let Q be its Goldie quotient ring . Let

**c(0) be the -set of regular elements of R . For a fixed positive integer , m, m < n , let g be the m- Gabriel filter defined as g = { all right ideals I in R / | R/I | < m} . Let cm = { all c ε R / | R/cRI | < m} and let c = { family of right ideals k of R / k п Cm ≠ ϕ } . Let xm = { all prime ideals p ε R / | R/p | = m } be the full set of m- prime ideals of R . Let for each p ε xm , c(p) have the usual meaning . We denote for each p ε xm the family gp = { family of right ideals J in g / J п c(p) ≠ ϕ } . Let $R^e$ = { q in Q / qJ < R , for some J in g } be the Gabriel extension ring of R at the Gabriel filter g. Then we have the following ;**

**(a) c ≤ g ≤ gp for each p ε xm .Hence g is a localising filter for cp for each p ε xm .**

**(b) g= c, hence cm is a right oreset .**

**(c) Moreover if Rc denotes the right quotient ring of R at cm then Rc equals the Gabriel quotient ring of R at the Gabriel filter g . .**

**(d) We also have that Xm ≤ maximal spec.Rc≤ prim (Rc)**

**Proof :- (a) First note that (a) is obvious .**

**(b) That (b) follows from (a) is also obvious. Denote by Rc the right quotient ring of R at cm .**

**(c) This is straightforward after observing that g=c implies that for every right ideal J in g we have that there is an element c of cm such that cε g . We leave the rest of the proof to the reader .**

**(d) This is also straightforward .**

**<u>(3.4) Corollary</u> :- Let Rc denote the right quotient ring of R at the right ore set cm then xm has the right intersection condition if every simple Rc module Rc / kc where kc is a maximal right ideal of Rc with r-ann(Rc/Kc) =qc is such that Rc/kc is a torsionfree Rc/qc module. In particular thus xm has the right intersection condition if R is a fully bounded noetherian ring .**

Proof :- We prove the if condition in two steps ..

(i) (step 1)Suppose under the hypothesis of the theorem xm does not have the right intersection condition . So assume k is a right ideal of R such that k п c(p) ≠ ф , for each p ε xm , . but k п cm = ф Choose without loss of generality k maximal such that k п cm =ф Then clearly R/k is a critical right R module with IR/kI ≥ m .

(ii) (step(2) Now consider the ring Rc and note that since k п cm = ф hence k Rc is a proper right ideal of Rc . Moreover k is a maximal right ideal of R such that kп cm = ф implies that kRc is a maximal such right ideal of Rc . Hence Rc/Kc is a simple right Rc –module Let qc = largest two sided ideal in k Rc . Hence by the hypothesis of the theorem Rc/kc is a torsion free right Rc/ qc module hence R/ k is also a torsion free right R/ q module Thus we must have that IR/kI = m = IR/q I Hence q ε xm . But from the beginning of step (1) we assumed that k п c(p) ≠ ф , for each p ε xm . Thus our assumption is contradicted. Hence we must have that the right intersection condition holds true .

Finally we end this paper by summing the main theorem of our paper that we have proved .

<u>(3.5) Main Theorem ( construction and definition ) :-</u> Let R be a prime noetherian ring with k.dim. (R) = n , n a positive integer. Let Q denote the Goldie quotient ring of R . For a fixed positive integer m, 0 <m < n , call xm = {all pεspec(R) / k-dim. (R/p))=m} a full set of prime ideals of R . Let cm = { all cε R / k-dim. (R/cR )< m } , and let g = {all right ideals J ≤ R / k-dim. (R/ J )< m } be the m-Gabriel filter in R . We then construct an extension ring R(m) of R at the m-Gabriel filter g with the following properties (i) R(m) is a right noetherian subring of Q with identity element 1 of R . (ii) If u(R(m)) is the set of units of R(m) then u(R(m) п R = cm precisely.

(iii) For a full set of m-prime ideals the set cm is a right ore set .

**References :-**